\documentclass[12pt]{article}
\usepackage{amsfonts}
\usepackage{amsmath}

\newcommand{\Q}{\mathbb{Q}}   
   
\newcommand{\Z}{\mathbb{Z}}   
\newcommand{\PR}{\mathbb{P}}

\newtheorem{thm}{Theorem}

\newtheorem{prop}[thm]{Proposition}

\newtheorem{lem}[thm]{Lemma}

\newcommand{\e}{\mathrm{e}}
\newcommand{\pf}{\noindent {\bf PROOF.} \quad}
\newcommand{\qed}{$\Box$}

\newcommand{\Li}{\mathrm{Li}} 
\newcommand{\ii}{\mathrm{i}} 
\newcommand{\Lf}{\mathrm{L}} 
\newenvironment{ack}{\centering{\bf{Acknowledgments}.}}{ } 

\begin{document}
\title{On certain combination of colored multizeta values}

\author{Matilde N. Lal\'{\i}n {\footnote{This material is based upon work supported by the National Science Foundation under agreement No. DMS-0111298.}}\\ \small{Institute for Advanced Study, School of Mathematics}\\ \small{1 Einstein Drive, Princeton, NJ 08540, USA}\\
\small{e-mail: \texttt{mlalin@math.ias.edu}}}

\date{}
\maketitle

\begin{abstract}
We reduce the sum $\sum_{k=1}^\infty \sum_{j=0}^{k-1} \frac{(-1)^{j+k+1}}{(2j+1)^mk}$ in terms of special values of the Dirichlet L-series in the character of conductor 4. This sum is a combination of colored zeta values. 
\end{abstract}

\section{Introduction}

The sums \[ \Li_{n,m}(1,1)=\sum_{0< j <k} \frac{1}{j^mk^n}\] (where $m,n$ are positive integers and $n$ is greater than 1 for convergence) were first studied by Euler \cite{E}, who found a closed formula by reducing them as rational combinations of products of values of the Riemann zeta function for the case when $m+n$ is odd. The simplest example is the identity
\[\zeta(3) = \sum_{0< j <k} \frac{1}{j k^2}\]
In general, multizeta values were studied in different contexts by several people, including Hoffman \cite{H}, Zagier \cite{Z}, Kontsevich \cite{K}, Broadhurst \cite{B}, Goncharov \cite{G1,G2,G3}, Drinfeld \cite{D}, and many others \footnote{The literature on the subject is very vast. We do not intend to do a survey. The interested reader is advised to follow the references in the aforementioned articles for more information.}. 

A natural generalization is to consider sums of the form{\footnote{We use the notation from multiple polylogarithms, although we will not study them in general.}} 
\[\Li_{n,m}(\zeta_l,\zeta_h)=\sum_{0< j <k} \frac{\zeta_h^j \zeta_l^k}{j^mk^n},\]
where $\zeta_h$, $\zeta_l$ are $l^{th}$- and $h^{th}$-roots of unity. We have the same conditions as before for $m$ and $n$, but $n$ is allowed to be equal to $1$ if $\zeta_l\not=1$, in which case the series converges but not absolutely. In such a situation, the sum is performed on the variable $j$ first.
   
These sums were first studied by Deligne \cite{De1,De2}, Goncharov \cite{G4,G2}, Racinet \cite{R,R2}, Bigotte, Jacob, Oussous, and Petitot \cite{BJOP} among others. Once again the reader is referred to the references in these works. 

One of the main problems concerning multizeta values and their generalizations is to understand and describe the relations among them. We will be concerned with explicit relations.

For the cases when $l=h=2$ and $m+n$ is odd, a summary of the results can be found in formula (75) of the work by Borwein, Bradley, and Broadhurst, \cite{BBB}, more precisely,
\begin{multline} \label{eq:bbb}
 \Li_{n,m}(\rho,\sigma) = \frac{1}{2} \left( - \Li_{m+n}(\rho \sigma) + (1+(-1)^n)\Li_n(\rho)\Li_m(\sigma) \right)\\ + \frac{(-1)^n}{2} \left( \binom{m+n-1}{n-1} \Li_{m+n}(\rho) + \binom{m+n-1}{m-1} \Li_{m+n} (\sigma) \right) \\
- \sum_{0 < k < \frac{m+n}{2} }\Li_{2k}(\rho \sigma) (-1)^n\left( \binom{m+n-2k-1}{n-1} \Li_{m+n-2k}(\rho) + \binom{m+n-2k-1}{m-1} \Li_{m+n-2k} (\sigma) \right),
\end{multline}
for $m+n$ odd, $\rho=\pm 1$, and $ \sigma=\pm 1$. 

Now, if we consider twisting by fourth roots of the unity, we note that
\[ \Li_{n,m}(\ii,\ii) - \Li_{n,m}(\ii,-\ii) + \Li_{n,m}(-\ii,\ii) - \Li_{n,m}(-\ii,-\ii) = 2^{2-n} \ii \sum_{0\leq j <k} \frac{(-1)^{j+k}}{(2j+1)^mk^n}.\]
When $n=1$ the sum is performed first in the variable $j$.

Here is our result:
\begin{thm} \label{bigresult} For $m \in \Z$, and $m$ odd, we have
\begin{multline}
 \sum_{k=1}^\infty \sum_{j=0}^{k-1} \frac{(-1)^{j+k+1}}{(2j+1)^mk} = m \Lf(\chi_{-4}, m+1)
 + \sum_{h=1}^\frac{m-1}{2} \frac{(-1)^h \pi^{2h}(2^{2h}-1)}{(2h)!} B_{2h} \Lf(\chi_{-4},m-2h+1). 
\end{multline}
where $\Lf(\chi_{-4},n)$ is the Dirichlet L-series in the odd character of conductor 4.
\end{thm}
Observe that it is possible to express $\Lf(\chi_{-4},n)$ in terms of special values of the polylogarithm $\Li_n$, namely, $\Lf(\chi_{-4},n) = \frac{\Li_n(\ii) - \Li_n(-\ii) }{2 \ii}$. Hence Theorem \ref{bigresult} becomes a result about reducing a particular combination of polylogarithms of depth 2 into a combination of polylogarithms of depth 1.

It should be also possible to reduce the more general sum $\sum_{0\leq j <k} \frac{(-1)^{j+k}}{(2j+1)^mk^n}$ for $m+n$ even and $n>1$ by using the same ideas that we are about to describe. However the computation becomes too complicated. 

\section{Idea of the proof} \label{idea}

The main idea of the proof was inspired by the following result by Murty and Sinha \cite{MS}.

\begin{thm}(4.1 in \cite{MS}) \label{ms} For $0<x<1$,
\begin{equation} \label{eq:ms}
-\frac{1}{x^s} +\zeta(s;x) = \sum_{r=0}^\infty \binom{-s}{r} \zeta(s+r) x^r,
\end{equation}
where $\zeta(s;x)= \sum_{n=0}^\infty \frac{1}{(n+x)^s}$ is the Hurwitz zeta function.
\end{thm}
This equality is proved by noting that
\[ \frac{1}{(n+x)^s} = \frac{1}{n^s} \sum_{r=0}^\infty \binom{-s}{r} \left( \frac{x}{n}\right)^r,\]  then summing over $n$, and interchanging the order of the sums.
Observe that this argument also works when $s$ is integral for $-1<x<0$.

For our proof, we first start by performing the change $j \rightarrow j-1$, and we obtain the following identity
\[ S:= \sum_{k=1}^\infty \sum_{j=0}^{k-1} \frac{(-1)^{j+k+1}}{(2j+1)^mk} = 
\sum_{k=1}^\infty \sum_{j=1}^{k}\frac{(-1)^{j+k}}{(2j-1)^mk}.\]
Now we may combine both sums in order to obtain
\[ 2S= \sum_{k=1}^\infty \sum_{j=0}^{k-1} \frac{(-1)^{j+k+1}}{(2j+1)^mk} - \sum_{k=1}^\infty \sum_{j=1}^{k} \frac{(-1)^{j+k+1}}{(2j-1)^mk}\]
\begin{equation} \label{eq:1}
=  \sum_{k=1}^\infty \frac{(-1)^{k+1}}{k} +\sum_{k=1}^\infty \frac{1}{(2k-1)^m k} +
 \sum_{k=1}^\infty \sum_{j=1}^{k-1}\frac{(-1)^{j+k+1}}{k} \left( \frac{1}{(2j+1)^m} - \frac{1}{(2j-1)^m}\right).
\end{equation}

The second term may be easily expressed as a combination of values of the Riemann zeta function and the alternating harmonic series:
\[ \sum_{k=1}^\infty \frac{1}{(2k-1)^m k} =2 \sum_{j=0}^{m-2}(-1)^j   \left( 1 - \frac{1}{2^{m-j}}\right) \zeta(m-j)  + (-1)^{m-1} \sum_{k=1}^\infty \left ( \frac{2}{2k-1} - \frac{1}{k}\right).\]

Thus we obtain
\[  \sum_{k=1}^\infty \frac{(-1)^{k+1}}{k} +\sum_{k=1}^\infty \frac{1}{(2k-1)^m k}= 3 \log 2 + 2 \sum_{j=0}^{m-2}(-1)^j   \left( 1 - \frac{1}{2^{m-j}}\right) \zeta(m-j).   \]

The third term in equation (\ref{eq:1}) is 
\begin{multline*}
\sum_{k=1}^\infty \sum_{j=1}^{k-1} \frac{(-1)^{j+k+1}}{k} \left( \frac{1}{(2j+1)^m} - \frac{1}{(2j-1)^m}\right)= \sum_{k=1}^\infty \sum_{j=1}^{k-1}\frac{(-1)^{j+k+1}}{k(2j)^m} \sum_{r=0}^\infty \binom{-m}{r} \frac{1-(-1)^r}{(2j)^r}\\ = - \sum_{r=1,\, r \, \mathrm{odd}}^\infty \binom{-m}{r} \frac{1}{2^{r+m-1}}\sum_{k=1}^\infty \sum_{j=1}^{k-1}\frac{(-1)^{j+k}}{kj^{r+m}}\\=  - \sum_{r=1,\, r \, \mathrm{odd}}^\infty \binom{-m}{r} \frac{1}{2^{r+m-1}} \Li_{1,r+m}(-1,-1)
\end{multline*}

We will use formula (\ref{eq:bbb}), which in this particular case implies
\begin{multline*} 
 \Li_{1,r+m}(-1,-1)=  \left( \frac{r+m}{2} -\frac{r+m+1}{2^{r+m+1}}\right) \zeta(r+m+1) -\zeta(r+m) \log 2\\- \sum_{k=1}^{\frac{r+m-2}{2}} \zeta(2k) \left( 1- \frac{1}{2^{r+m-2k}}\right) \zeta(r+m+1-2k).
\end{multline*}
At this point it is very clear that we will have to work with equations that are similar to (\ref{eq:ms}), but not exactly the same.

\section{Some helpful results}
We have the following extension of Theorem \ref{ms}.
\begin{prop}  \label{prop:1} For $-1<x<1$ but $x\not =0$ and $s, t$ positive integers,
 \begin{multline} \label{eq:p1}
 \sum_{r=0}^\infty \binom{-s}{r} \zeta(r+s+t) x^{r}  = \sum_{h=2}^t \frac{(-1)^{t-h}}{x^{s+t-h}}\binom{s+t-h-1}{t-h} \zeta(h) \\
 + \frac{(-1)^{t-1}}{x^{s+t-1}} \binom{s+t-2}{t-1} \sum_{n=1}^\infty \left( \frac{1}{n} - \frac{1}{n+x} \right)\\ + (-1)^t \sum_{h=2}^s \frac{1}{x^{s+t-h}} \binom{s+t-h-1}{t-1} \sum_{n=1}^\infty \frac{1}{(n+x)^h}.
\end{multline} 

\end{prop}

\pf The key observation is the same as in Theorem  \ref{ms}.
\[  \sum_{n=1}^\infty \frac{1}{n^t (n+x)^s} = \sum_{n=1}^\infty \frac{1}{n^{t+s}}\sum_{r=0}^\infty \binom{-s}{r} \left( \frac{x}{n}\right)^r  =  \sum_{r=0}^\infty \binom{-s}{r} \zeta(s+t+r) x^r.\]

Now we will need to use
 \begin{equation} \label{eq:2}
 \frac{1}{n^t (n+x)^s} = \sum_{h=1}^t \frac{(-1)^{t-h}}{n^h x^{s+t-h}} \binom{s+t-h-1}{t-h} + (-1)^t \sum_{h=1}^s \frac{x^{h-s-t}}{(n+x)^h} \binom{s+t-h-1}{t-1},
 \end{equation}
which may be proved by induction on $t$. It is true for $t=1$ because
\[\frac{1}{n (n+x)^s} =\frac{1}{nx^s} - \sum_{h=1}^s \frac{x^{h-s-1}}{(n+x)^h} \]
Assume we know it for $t$, then
\[\frac{1}{n^{t+1} (n+x)^s} = \sum_{h=1}^t \frac{(-1)^{t-h}}{n^{h+1} x^{s+t-h}} \binom{s+t-h-1}{t-h} + (-1)^t \sum_{h=1}^s \frac{x^{h-s-t}}{n(n+x)^h} \binom{s+t-h-1}{t-1}\]
 \[=  \sum_{h=2}^{t+1} \frac{(-1)^{t+1-h}}{n^{h} x^{s+t+1-h}} \binom{s+t-h}{t+1-h} + (-1)^t \sum_{h=1}^s x^{h-s-t}  \left(\frac{1}{nx^h} - \sum_{l=1}^h \frac{x^{l-h-1}}{(n+x)^l}  \right) \binom{s+t-h-1}{t-1}\]
\[=  \sum_{h=1}^{t+1} \frac{(-1)^{t+1-h}}{n^{h} x^{s+t+1-h}} \binom{s+t-h}{t+1-h} + (-1)^{t+1} \sum_{l=1}^s \frac{x^{l-s-t-1}}{(n+x)^l}   \sum_{h=l}^s \binom{s+t-h-1}{t-1}\]
This proves equation (\ref{eq:2}) and from there is very easy to establish the Proposition. \qed

\begin{prop} \label{prop:2} For $2k+1 \geq s$,
\begin{multline}  
  \sum_{r=0}^\infty \binom{-s}{r+2k+1-s} \zeta(r+2) x^{r}\\
  = \sum_{h=0}^{s-2} (-1)^{s-h} \binom{2k-1}{h} (\zeta(s-h;x) x^{s-h-2} - x^{-2}) \\
  + \binom{2k-1}{s-1}\frac{1}{x}\sum_{n=1}^\infty \left( \frac{1}{n} - \frac{1}{n+x}\right).
  \end{multline}
\end{prop}
\pf We start by integrating  identity (\ref{eq:ms}) with $s=2$ between 0 and $x$. Then we multiply by $x^{2k-1}$:
\[x^{2k-1}\sum_{n=1}^\infty \left( \frac{1}{n} - \frac{1}{n+x}\right)  = \sum_{r=0}^\infty \binom{-2}{r} \zeta(r+2)\frac{x^{r+2k}}{r+1}.\]
Now we differentiate $s-1$ times,
\begin{multline*}
(s-1)! \sum_{h=0}^{s-2} (-1)^{s-h} \binom{2k-1}{h} (\zeta(s-h;x) x^{2k-h-1} - x^{2k-1-s})\\
 + (s-1)!\binom{2k-1}{s-1} x^{2k-s} \sum_{n=1}^\infty \left( \frac{1}{n} -\frac{1}{n+x}\right)\\=  (s-1)! \sum_{r=0}^\infty \binom{-s}{r+2k+1-s} \zeta(r+2) x^{r+2k+1-s},
\end{multline*}
we obtain the result by dividing by $x^{2k-s+1}$ and $(s-1)!$. \qed

In the proof of Theorem \ref{bigresult} we are going to need to simplify certain series involving even values of the Riemann zeta function. The following Lemma will be specially useful.
\begin{lem} \label{lemma} If $h>0$,
\begin{equation}
\sum_{k=1}^\infty \frac{\zeta(2k)}{2^{2k-1}} \binom{2k-1}{h} =-  \frac{(\ii \pi)^{h+1}(2^{h+1}-1)}{(h+1)!} B_{h+1} + (-1)^h
\end{equation}
where the $B_n$ are the Bernoulli numbers given by $\frac{t}{\e^t-1} = \sum_{n=0}^\infty \frac{B_n t^n}{n!}$.

If $h=0$,
\begin{equation}
\sum_{k=1}^\infty \frac{\zeta(2k)}{2^{2k-1}} = 1. 
\end{equation}

\end{lem}
\pf
First recall that
\[ \zeta(2k) =\frac{(-1)^{k-1}B_{2k}(2\pi)^{2k}}{2 (2k)!}.\]

We need to compute
\begin{equation}\label{eq:bernoulli}
\sum_{k=1}^\infty \frac{\zeta(2k)}{2^{2k-1}} \binom{2k-1}{h} = - \sum_{k=1}^\infty \frac{B_{2k}(\ii\pi)^{2k}}{(2k)!} \binom{2k-1}{h}.
\end{equation}
If $h>0$ we may write,
\[= - \sum_{n=h+1}^\infty \frac{B_{n}(\ii\pi)^{n}}{n!} \binom{n-1}{h} \]
\[ = -\frac{(\ii \pi)^{h+1}}{h!}  \sum_{n=h+1}^\infty \frac{B_{n}(\ii\pi)^{n-h-1}}{n!} (n-1) \dots (n-h)\]
\[ = -\frac{(\ii \pi)^{h+1}}{h!} \frac{\partial^{h}}{\partial t^{h}} \left . \left( \frac{1}{\e^t-1} -\frac{1}{t} \right) \right|_{t=\ii \pi}\] 
\[=  -\frac{(\ii \pi)^{h+1}}{h!} \frac{\partial^{h}}{\partial t^{h}} \left . \left( \frac{1}{\e^{t+\ii \pi}-1}\right) \right|_{t=0} +(-1)^h \] 
\[=  \frac{(\ii \pi)^{h+1}}{h!} \frac{\partial^{h}}{\partial t^{h}} \left . \left( \frac{1}{\e^{t}+1}\right) \right|_{t=0} +(-1)^h \] 
\[=  \frac{(\ii \pi)^{h+1}}{h! 2} E_h(0) + (-1)^h \]
where the $E_n(x)$ are the Euler polynomials given by $\frac{2\e^{xt}}{e^t+1} = \sum_{n=0}^\infty \frac{E_n(x)t^n}{n!}$. Now we use that 
\[E_{n}(0) = - \frac{2(2^{n+1}-1) B_{n+1}}{n+1},\]
(see page 805 of \cite{AS}). Then equation (\ref{eq:bernoulli}) becomes
\[= -  \frac{(\ii \pi)^{h+1}(2^{h+1}-1)}{(h+1)!} B_{h+1} + (-1)^h.\] 

For the case when $h=0$, equation (\ref{eq:bernoulli}) becomes
\begin{multline}
\sum_{k=1}^\infty \frac{\zeta(2k)}{2^{2k-1}}  = - \sum_{n=0}^\infty \frac{B_{n}(\ii\pi)^{n}}{n!} + B_{0} + B_{1}\ii\pi= -\frac{\ii\pi}{\e^{\ii\pi}-1} + 1 -\frac{\ii \pi}{2} =1.
\end{multline}
\qed

\section{The conclusion of the proof}

We will now proceed to finish the proof of Theorem \ref{bigresult}. Recall that from section \ref{idea} we know that
\begin{multline}\label{eq:3}
2S= 3 \log 2 + 2 \sum_{j=0}^{m-2}(-1)^j   \left( 1 - \frac{1}{2^{m-j}}\right) \zeta(m-j) \\ 
-\sum_{r=1,\, r \, \mathrm{odd}}^\infty \binom{-m}{r} \frac{1}{2^{r+m-1}} \left(\left( \frac{r+m}{2} -\frac{r+m+1}{2^{r+m+1}}\right) \zeta(r+m+1) -\zeta(r+m) \log 2 \right.\\\left.- \sum_{k=1}^{\frac{r+m-2}{2}} \zeta(2k) \left( 1- \frac{1}{2^{r+m-2k}}\right) \zeta(r+m+1-2k) \right). 
\end{multline}
We proceed to evaluate the terms in the second and third lines of equation (\ref{eq:3}). First,
\[- \sum_{r=1,\, r \, \mathrm{odd}}^\infty \binom{-m}{r} \frac{r+m}{2^{r+m}} \zeta(r+m+1) = - \frac{m}{2^{m+1}} \sum_{r=0}^\infty \binom{-(m+1)}{r}  \zeta(r+m+1)\frac{1-(-1)^r}{2^{r}}= m, \]
by Theorem \ref{ms}.

The next term is
\[\sum_{r=1,\, r \, \mathrm{odd}}^\infty \binom{-m}{r}\frac{r+m+1}{2^{2r+2m}} \zeta(r+m+1)\]
\[=\sum_{r=1,\, r \, \mathrm{odd}}^\infty \binom{-m}{r}\frac{r+m}{2^{2r+2m}} \zeta(r+m+1) + \sum_{r=1,\, r \, \mathrm{odd}}^\infty \binom{-m}{r}\frac{1}{2^{2r+2m}} \zeta(r+m+1)\]
\[ = \frac{m}{2^{2m+1}} \sum_{r=0}^\infty \binom{-(m+1)}{r}  \zeta(r+m+1)\frac{1-(-1)^r}{2^{2r}} +\frac{1}{2^{2m+1}} \sum_{r=0}^\infty \binom{-m}{r}\zeta(r+m+1)\frac{1-(-1)^r}{2^{2r}}\]
by Theorem \ref{ms} and Proposition \ref{prop:1}:
\[= 2m \Lf(\chi_{-4}, m+1) - 2m + 2-3\log 2 + 2m-2 - 2\sum_{j=1}^{\frac{m-1}{2}} \left( \Lf(\chi_{-4},2j) + \left( 1- \frac{1}{2^{2j+1}}\right) \zeta(2j+1) \right)\]
\[=2m \Lf(\chi_{-4}, m+1) -3\log 2 - 2\sum_{j=1}^{\frac{m-1}{2}} \left( \Lf(\chi_{-4},2j) + \left( 1- \frac{1}{2^{2j+1}}\right) \zeta(2j+1) \right).\]

The third term in the second line of equation (\ref{eq:3}) corresponds to
\[\log 2 \sum _{r=1,\, r \, \mathrm{odd}}^\infty \binom{-m}{r} \zeta(r+m) \frac{1}{2^{r+m-1}} =\frac{\log 2}{2^m} \sum _{r=0}^\infty \binom{-m}{r} \zeta(r+m) \frac{1-(-1)^r}{2^r} = -\log 2.\]

Finally, each of the terms in the last line of equation (\ref{eq:3}) equals
\[ \zeta(2k) \sum_{r=2k+2-m,\, r \, \mathrm{odd}}^\infty \binom{-m}{r}\left( \frac{1}{2^{r+m-1}} -\frac{1}{2^{2r+2m-2k-1}}\right) \zeta(r+m+1-2k) \]
set $t=r-2k+m-1$,
\[= \zeta(2k) \sum_{t=1,\, t \, \mathrm{odd}}^\infty \binom{-m}{t+2k-m+1}\left( \frac{1}{2^{t+2k}} -\frac{1}{2^{2t+2k+1}}\right) \zeta(t+2) \]
\[= \frac{ \zeta(2k)}{2^{2k+1}} \sum_{t=0}^\infty \binom{-m}{t+2k-m+1}  \zeta(t+2) \frac{1-(-1)^t}{2^{t}} \]
\[ -  \frac{ \zeta(2k)}{2^{2k+2}} \sum_{t=0}^\infty \binom{-m}{t+2k-m+1}  \zeta(t+2) \frac{1-(-1)^t}{2^{2t}}.  \]
By Proposition \ref{prop:2},
\[= -\frac{\zeta(2k)}{2^{2k-1}} \binom{2k-1}{m-1} (2 \log 2-1) - \frac{\zeta(2k)}{2^{2k-1}} \sum_{j=0}^\frac{m-3}{2} \binom{2k-1}{2j+1}\]
\[- \frac{\zeta(2k)}{2^{2k-1}} \sum_{j=0}^\frac{m-3}{2} \binom{2k-1}{2j} \left( 2 \left(1- \frac{1}{2^{m-2j}}\right) \zeta(m-2j) -1 \right)\]
\[ + \frac{\zeta(2k)}{2^{2k-1}} \binom{2k-1}{m-1} (3 \log 2-2) \]
\[+ \frac{\zeta(2k)}{2^{2k-2}} \sum_{j=0}^\frac{m-3}{2} \binom{2k-1}{2j} \left( \left(1- \frac{1}{2^{m-2j}}\right) \zeta(m-2j) -1 \right)\]
\[-  \frac{\zeta(2k)}{2^{2k-2}} \sum_{j=0}^\frac{m-3}{2} \binom{2k-1}{2j+1} \left(\Lf(\chi_{-4},m-2j-1) -1 \right)\]
\[= \frac{\zeta(2k)}{2^{2k-1}} \binom{2k-1}{m-1}\log 2 - \frac{\zeta(2k)}{2^{2k-1}} \sum_{j=0}^{m-1} (-1)^j\binom{2k-1}{j}\]
\[-  \frac{\zeta(2k)}{2^{2k-2}} \sum_{j=0}^\frac{m-3}{2} \binom{2k-1}{2j+1} \Lf(\chi_{-4},m-2j-1).\]

At this point it will be necessary to use Lemma \ref{lemma}. The last line in equation (\ref{eq:3}) is
\[\sum_{k=1}^\infty \left( \frac{\zeta(2k)}{2^{2k-1}} \binom{2k-1}{m-1}\log 2 - \frac{\zeta(2k)}{2^{2k-1}} \sum_{j=0}^{m-1} (-1)^j\binom{2k-1}{j} \right .\]
\[\left .-  \frac{\zeta(2k)}{2^{2k-2}} \sum_{j=0}^\frac{m-3}{2} \binom{2k-1}{2j+1} \Lf(\chi_{-4},m-2j-1)  \right)\]
\[=  \log 2 -1 - \sum_{j=1}^{m-1} \left( \frac{(-\ii \pi)^{j+1} (2^{j+1}-1)}{(j+1)!} B_{j+1} + 1 \right)\] 
\[+2 \sum_{j=0}^\frac{m-3}{2} \left( \frac{(\ii \pi)^{2j+2} (2^{2j+2}-1)}{(2j+2)!} B_{2j+2} + 1 \right ) \Lf(\chi_{-4},m-2j-1).  \]
Observe that the above equation is still true for $m=1$.

Setting together all the terms for equation \ref{eq:3}, we obtain, 
\[ 2S =m +2m \Lf(\chi_{-4}, m+1) - 2\sum_{j=1}^{\frac{m-1}{2}} \left( \Lf(\chi_{-4},2j) + \left( 1- \frac{1}{2^{2j}}\right) \zeta(2j) \right) -\log2 \]
\[+\log 2  -1 -  \sum_{j=1}^{m-1}\left(   \frac{(-\ii \pi)^{j+1}(2^{j+1}-1)}{(j+1)!} B_{j+1} + 1\right)\]
\[ + 2\sum_{j=0}^\frac{m-3}{2}\left(  \frac{(-1)^{j+1} \pi^{2j+2}(2^{2j+2}-1)}{(2j+2)!} B_{2j+2} +1 \right) \Lf(\chi_{-4},m-2j-1).\] 
Finally,
\[ 2 S= 2m \Lf(\chi_{-4}, m+1) +2 \sum_{h=1}^\frac{m-1}{2} \frac{(-1)^h \pi^{2h}(2^{2h}-1)}{(2h)!} B_{2h} \Lf(\chi_{-4},m-2h+1).\] 
\qed

\section{An application to Mahler measure}

A motivation for studying this particular sum comes from the world of Mahler measure (see for instance, \cite{B1}). Some formulas for Mahler measure of multivariate polynomials are computed in  \cite{L, L2}. These formulas express the Mahler measure of certain families of polynomials in terms of special values of Dirichlet L-functions, the Riemann zeta function, and sometimes sums of the kind  $\sum_{0\leq j <k} \frac{(-1)^{j+k+1}}{(2j+1)^3k^{2l+1}} $.  However, one would expect to obtain formulas which would only depend on polylogarithms of depth one. Theorem \ref{bigresult} allows us to simplify one of these Mahler measure formulas.

\begin{thm}
\begin{equation}
\pi^3 m\left(1+x+\left(\frac{1-x_1}{1+x_1}\right) (1+y) z\right) = 24 \Lf(\chi_{-4}, 4). 
\end{equation}
\end{thm}

\pf From \cite{L} we know that
\[\pi^3 m\left(1+x+\left(\frac{1-x_1}{1+x_1}\right) (1+y) z\right) = 2 \pi^2 \Lf(\chi_{-4},2) + 8 \sum_{k=1}^\infty \sum_{j=0}^{k-1} \frac{(-1)^{j+k+1}}{(2j+1)^3k}.\]
Applying Theorem \ref{bigresult} we obtain the statement. \qed

Results for the general case of $\sum_{0\leq j <k} \frac{(-1)^{j+k}}{(2j+1)^mk^n}$ would lead to the simplification of all the formulas in \cite{L2}.

\bigskip
\begin{ack}
The author wishes to thank M. Ram Murty for helpful discussions and for providing a copy of \cite{MS}.
\end{ack}

\end{document}